\documentstyle[12pt]{article}

\expandafter\chardef\csname pre amssym.def at\endcsname=\the\catcode`\@
\catcode`\@=11

\def\undefine#1{\let#1\undefined}
\def\newsymbol#1#2#3#4#5{\let\next@\relax
 \ifnum#2=\@ne\let\next@\msafam@\else
 \ifnum#2=\tw@\let\next@\msbfam@\fi\fi
 \mathchardef#1="#3\next@#4#5}
\def\mathhexbox@#1#2#3{\relax
 \ifmmode\mathpalette{}{\m@th\mathchar"#1#2#3}%
 \else\leavevmode\hbox{$\m@th\mathchar"#1#2#3$}\fi}
\def\hexnumber@#1{\ifcase#1 0\or 1\or 2\or 3\or 4\or 5\or 6\or 7\or 8\or
 9\or A\or B\or C\or D\or E\or F\fi}

\font\tenmsa=msam10 scaled 1000
\font\sevenmsa=msam7 scaled 1000
\font\fivemsa=msam5 scaled 1000
\newfam\msafam
\textfont\msafam=\tenmsa
\scriptfont\msafam=\sevenmsa
\scriptscriptfont\msafam=\fivemsa
\edef\msafam@{\hexnumber@\msafam}
\mathchardef\dabar@"0\msafam@39
\def\dashrightarrow{\mathrel{\dabar@\dabar@\mathchar"0\msafam@4B}}
\def\dashleftarrow{\mathrel{\mathchar"0\msafam@4C\dabar@\dabar@}}

\def\ulcorner{\delimiter"4\msafam@70\msafam@70 }
\def\urcorner{\delimiter"5\msafam@71\msafam@71 }
\def\llcorner{\delimiter"4\msafam@78\msafam@78 }
\def\lrcorner{\delimiter"5\msafam@79\msafam@79 }
\def\yen{{\mathhexbox@\msafam@55 }}
\def\checkmark{{\mathhexbox@\msafam@58 }}
\def\circledR{{\mathhexbox@\msafam@72 }}
\def\maltese{{\mathhexbox@\msafam@7A }}

\font\tenmsb=msbm10 scaled 1000
\font\sevenmsb=msbm7 scaled 1000
\font\fivemsb=msbm5 scaled 1000
\newfam\msbfam
\textfont\msbfam=\tenmsb
\scriptfont\msbfam=\sevenmsb
\scriptscriptfont\msbfam=\fivemsb
\edef\msbfam@{\hexnumber@\msbfam}
\def\Bbb#1{\fam\msbfam\relax#1}
\def\widehat#1{\setboxz@h{$\m@th#1$}%
 \ifdim\wdz@>\tw@ em\mathaccent"0\msbfam@5B{#1}%
 \else\mathaccent"0362{#1}\fi}
\def\widetilde#1{\setboxz@h{$\m@th#1$}%
 \ifdim\wdz@>\tw@ em\mathaccent"0\msbfam@5D{#1}%
 \else\mathaccent"0365{#1}\fi}

\font\teneufm=eufm10 scaled 1000 
\font\seveneufm=eufm7 scaled 1000
\font\fiveeufm=eufm5 scaled 1000
\newfam\eufmfam
\textfont\eufmfam=\teneufm
\scriptfont\eufmfam=\seveneufm
\scriptscriptfont\eufmfam=\fiveeufm

\catcode`\@=\csname pre amssym.def at\endcsname

\def\ylth{0.5pt}

\def\yttv#1#2#3#4#5#6#7#8#9{%
  \smash{%
    \rule[1.5ex]{#1ex}{\ylth}\kern-#1ex%
    \rule[0.5ex]{#1ex}{\ylth}\kern-#1ex%
    \rule[-0.5ex]{#2ex}{\ylth}\kern-#2ex%
    \rule[-1.5ex]{#3ex}{\ylth}\kern-#3ex%
    \rule[-2.5ex]{#4ex}{\ylth}\kern-#4ex}%
  \raisebox{1.5ex}{%
    \rule[-#5ex]{\ylth}{#5ex}\kern-\ylth\kern1ex%
    \rule[-#5ex]{\ylth}{#5ex}\kern-\ylth\kern1ex%
    \rule[-#6ex]{\ylth}{#6ex}\kern-\ylth\kern1ex%
    \rule[-#7ex]{\ylth}{#7ex}\kern-\ylth\kern1ex%
    \rule[-#8ex]{\ylth}{#8ex}\kern-\ylth\kern1ex%
    \rule[-#9ex]{\ylth}{#9ex}\kern-\ylth\kern-5ex}%
  \kern#1ex\rule[1.5ex]{\ylth}{\ylth}\kern-\ylth}
\def\ytt#1#2#3#4#5#6#7#8#9{~\yttv#1#2#3#4#5#6#7#8#9~}

\newcommand{\mc}[1]{{\cal #1}}
\newcommand{\A}{{\cal A}}
\newcommand{\otA}{\otimes _{\!{\tiny \mbox{$\A$}}}}
\newcommand{\br}{{\bf r}}
\newcommand{\cont}{{^{\rm c}}}
\newcommand{\contr}{{\rm c}}
\newcommand{\Rda}{\hat{R}}
\newcommand{\Rdam}{{\hat{R}^{-1}}}
\newcommand{\Rup}{{\check{R}}}
\newcommand{\Rupm}{{\check{R}^{-}}}
\newcommand{\Rlim}{{\grave{R}^{-}}}
\newcommand{\Mor}{{\rm Mor}}
\newcommand{\im}{{\rm im}}

\newcommand{\Gamm}{\mbox{$\mit{\Gamma}$}}
\def\ot{\otimes}
\def\ott{{\otimes}}
\def\C{{\Bbb C}}
\newcommand{\QM}{q^{-1}}
\newcommand{\RM}{r^{-1}}
\newcommand{\lam}{{\lambda}}
\def\eps{\epsilon}

\newcommand{\fettc}{\contr}

\newcommand{\gten}{{\Gamm^{\otimes}}}
\newcommand{\gd}{{\Gamm^{\land}}}
\def\dota{\dot{A_k}}

\def\da{\overline{A}_k}
\def\Rchm{\check{R}^{-}}

\def\zchk{\check{z}_{(k)}}
\def\pch{\check{\pi}_\mu}

\begin{document}
\begin{center}
{\LARGE 
Exterior Algebras Related to the Quantum Group $\mc{O}(O_q(3))$
}\\
\vskip 3em
{\large Istv{\'a}n  Heckenberger\footnotemark[1]\ \quad \quad \quad
Axel Sch{\"u}ler\footnotemark[2]}\\
\vskip 1.5em
{\large Universit{\"a}t Leipzig,
Mathematisches Institut,
Augustusplatz 9-11,
04109 Leipzig, Germany}
\end{center}

\footnotetext[1]{e-mail: heckenbe@mathematik.uni-leipzig.de}
\footnotetext[2]{e-mail: schueler@mathematik.uni-leipzig.de\\
Supported by the Deutsche Forschungsgemeinschaft}

\begin{abstract}
For the 9-dimensional bicovariant differential calculi
on the quantum group $\mc{O}(O_q(3))$ several kinds
of exterior algebras are examined. The corresponding
dimensions, bicovariant subbimodules and eigenvalues of the
antisymmetrizer are given. Exactly one of the
exterior algebras studied by the authors has a unique left
invariant form with maximal degree.
\end{abstract}

\section{Introduction}

A general framework of bicovariant differential calculus on quantum groups
was given by the pioneering work of 
Woronowicz \cite{a-Woro2}. Bicovariant first order differential calculi on 
$q$-deformed simple Lie groups were constructed, studied, and classified by 
many authors, see \cite{a-Woro2,a-CSchWW1, J, SS}.
Recently several problems connected 
with 
higher order differential forms (exterior algebras) were studied. But only in 
case of the simplest examples like  $GL_q(N)$ and $SL_q(N)$ the structures
of the left-invariant and of the bi-invariant exterior algebras are known, 
see \cite{H,a-Schueler,T}.
For the orthogonal and  symplectic quantum groups even the
existence 
of an analogue of the volume form, i.\,e. a single form of maximal degree, was 
open. There are two  purposes of this paper. Firstly we show that for the 
$N^2$-dimensional 
bicovariant first order differential calculi on $O_q(N)$ Woronowicz' 
left-invariant external algebra is infinite dimensional. In other words, for 
each positive integer $k$ there exists a nonzero $k$-form. Secondly, in case
of 
$O_q(3)$ we discuss alternative constructions of exterior algebras which
yield a {\em finite} differential complex. For several choices of the ideal
of symmetric forms in $\gten$ we compute the dimensions of spaces of
left-invariant $k$-forms. Exactly one of those exterior algebras has a unique
left invariant form of maximal degree.

Let $\A$ be the Hopf algebra $\mc{O}(O_q(N))$ as defined in
\cite{a-FadResTak1}.
The fundamental matrix corepresentation of $\A$ is denoted by $u$.
We use the symbol $\Delta$ for the comultiplication and
Sweedler's notation $\Delta(a)=\sum a_{(1)}\ot a_{(2)}$.
The two-sided ideal of an algebra generated by a set $\{a_i\,|\,i\in I\}$
is denoted by $\langle a_i\,|\, i\in I\rangle$.
Let $v$  be a corepresentation  of
$\A$. As usual $v^\fettc$
denotes the contragredient corepresentation  of $v$. The
space of intertwiners of corepresentations $v$ and $w$ is  $\Mor(v,w)$.
We write $\Mor(v)$
for $\Mor(v,v)$.
Lower indices of a matrix $A$ always refer
to the components of a tensor product where $A$ acts (`leg numbering').
The unit matrix is
denoted by $I$.
As usual $\Rda$ and $C$ stand for the corresponding $\Rda$-matrix and the
metric, see \cite{a-FadResTak1}.

\section{Exterior Algebras}

Let $\A$ be a coquasitriangular Hopf algebra \cite{b-KS} with universal
$r$-form $\br$ and $v=(v^i_j)$ an arbitrary $n$-dimensional
corepresentation of $\A$. Let $\Gamm :=\Gamm (v)$ be the
corresponding bicovariant bimodule such that
$\{\theta_{ij}\,|\,i,j=1,\ldots,n\}$ is a basis of the vector space
of left invariant 1-forms and
\begin{eqnarray*}
\Delta _R(\theta _{ij})&\!\!=&\!\!\sum_{k,l=1}^n
\theta _{kl}\otimes v^k_i(v\cont)^l_j,\\
\theta _{ij}a&\!\!=&\!\!
\sum_{k,l=1}^n a_{(1)}\br(v^k_i,a_{(2)})\br(a_{(3)},v^j_l)\theta _{kl}.
\end{eqnarray*}
Let $\Gamm^{\otimes k}$ denote the $k$-fold algebraic tensor product
$\Gamm\otA\Gamm\otA\cdots\otA\Gamm$ ($k$ factors) of $\Gamm$
and $\Gamm^\otimes :=\sum _{k=0}^\infty \Gamm^{\otimes k}$.
Let $\mc{S}$ be a bicovariant subbimodule and two-sided ideal of
$\Gamm ^\otimes$. Then
$\Gamm ^\wedge :=\Gamm ^\otimes/\mc{S}$ is called an
\textit{exterior algebra of} $\Gamm$.
Since $\Gamm ^\otimes$ has a ${\Bbb Z}$-gradation, we require
$\mc{S}=\bigoplus_{k=0}^\infty\mc{S}^k$.
(Traditionally, the space $\mc{S}$ of symmetric forms has some additional
properties but we don't want to consider them at the moment.)
The general theory of bicovariant bimodules gives
$\Gamm^{\otimes k}=\A\Gamm_L^{\otimes k}$ and
$\mc{S}^k=\A\mc{S}_L^k$, where
$\Gamm_L^{\otimes k}=\{\rho\in\Gamm^{\otimes k}\,|\,
 \Delta_L(\rho)=1\otimes \rho\}$ and
$\mc{S}_L^k =\{\rho\in\mc{S}^k\,|\,\Delta_L(\rho)=1\otimes \rho\}$
are the corresponding left-invariant subspaces.
In the nontrivial cases we have $\mc{S}_L^0=\mc{S}_L^1=\{0\}$.

Let $\sigma$ be the canonical braiding of $\Gamm\otA \Gamm$.
Recall that $\sigma $ is a homomorphism of bicovariant bimodules.
Let $\Rlim=(\Rlim{}^{ij}_{kl})$ be the complex matrix with entries
$\Rlim{}^{ij}_{kl}:=\br(v^j_k,v\cont{}^i_l)$ and let $\Rlim_i$ be this matrix
acting on the $i$-th and $(i+1)$-th component of a vector.
Let us define a $k$-twist $b_k$ by the formula
\begin{equation}\label{eq-bk}
b_k :=(\Rlim_2\Rlim_4\cdots\Rlim_{2k-2})(\Rlim_3\cdots\Rlim_{2k-3})\cdots\\
\cdots (\Rlim_{k-1}\Rlim_{k+1})\Rlim_k.
\end{equation}
As an example we draw the corresponding pictures for $\Rlim$ and $b_4$
in the language of braids.
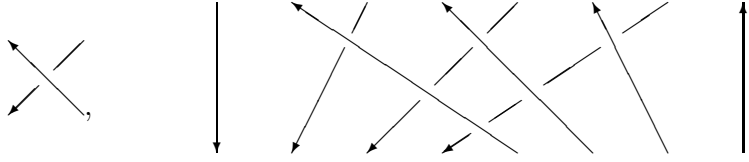
\begin{displaymath}
\setlength{\unitlength}{1cm}
\begin{picture}(1,1.5)
\put(1,0.5){\vector(-1,1){1}}
\put(.4,.9){\vector(-1,-1){.4}}
\put(.6,1.1){\line(1,1){.4}}
\end{picture}
\raisebox{.5cm}{,}\qquad\qquad
\begin{picture}(7,2)
\put(4,0){\vector(-3,2){3}}
\put(5,0){\vector(-1,1){2}}
\put(6,0){\vector(-1,2){1}}
\put(7,0){\vector(0,2){2}}
\put(0,2){\vector(0,-1){2}}
\put(2,2){\line(-1,-2){.2}}
\put(1.7,1.4){\vector(-1,-2){.7}}
\put(4,2){\line(-1,-1){.4}}
\put(3.4,1.4){\line(-1,-1){.5}}
\put(2.7,.7){\vector(-1,-1){.7}}
\put(6,2){\line(-3,-2){.6}}
\put(5.1,1.4){\line(-3,-2){.75}}
\put(4.05,.7){\line(-3,-2){.43}}
\put(3.38,.25333){\vector(-3,-2){.38}}
\end{picture}
\end{displaymath}
Since $\Rlim\in\Mor(v\otimes v\cont,v\cont \otimes v)$, we conclude that
$b_k\in\Mor(v^{\otimes k} v\cont{}^{\otimes k},(v v\cont)^{\otimes k})$.
Moreover, both $\Rlim$ and $b_k$ are invertible.

Let $P$ be a projection onto a subspace of $\Gamm^{\otimes k}$.
This subspace is $\Delta_R$-invariant iff $P\in\Mor((vv\cont)^{\otimes k})$.
For an endomorphism $T$ of the bicovariant bimodule $\Gamm^{\otimes k}$
let
\begin{equation}\label{eq-Tpunkt}
\dot{T}:=b_k^{-1}Tb_k.
\end{equation}
If there are projections
$P'\in\Mor(v^{\otimes k})$ and $P''\in\Mor((v\cont)^{\otimes k})$ such that
$\dot{P}=P'\odot P''$ then $\A\cdot P(\Gamm^{\otimes k}_L)$ is a
bicovariant subbimodule of $\Gamm^{\otimes k}$ where
\begin{displaymath}
(P'\odot P'')^{i_1,\ldots,i_{2k}}_{j_1,\ldots,j_{2k}}=
(P')^{i_1,\ldots,i_{k}}_{j_1,\ldots,j_{k}}
(P'')^{i_{k+1},\ldots,i_{2k}}_{j_{k+1},\ldots,j_{2k}}.
\end{displaymath}

This method gives a large class (in some cases all) of bicovariant
subbimodules of $\Gamm^{\otimes k}$ and the corresponding bimodule
homomorphisms. {}From now on we suppose that
\begin{displaymath}
\mc{S}^k=\A b_k(\sum_iP'_i\odot P''_i)b_k^{-1}(\Gamm^{\otimes k}_L),
\end{displaymath}
where $P'_i\in\Mor(v^{\otimes k})$ and $P''_i\in\Mor(v\cont{}^{\otimes k})$
are arbitrary projections.

\section{Exterior Algebras for $\mc{O}(SL_q(N))$}

Let us recall some results of \cite{a-Schueler}.
If $\A=\mc{O}(SL_q(N))$ and $v=u$ is the fundamental corepresentation,
then there are two canonical methods to
construct exterior algebras fitting in the above concept:

$a)$ $\mc{S}=\langle \ker(I-\sigma)\rangle $,

$b)$ $\mc{S}=\langle \ker A_k\,|\, k\geq 2\rangle $.

Here $A_k$ denotes the $k$-th antisymmetrizer on $\Gamm^{\otimes k}$
defined in \cite{a-Woro2}.

Both definitions give the same algebra $\Gamm^\wedge$.
The dimension of the vector space $\Gamm_L^{\wedge k}$ is
${N^2 \choose k}$, hence there is a unique left
invariant form with maximal degree.
The elements of $\Gamm^{\wedge k}$ can be given by
\begin{displaymath}
\A b_k(\sum_iP_i\odot \tilde{P}_i)b_k^{-1}(\Gamm^{\otimes k}_L),
\end{displaymath}
where the projections $(P_i\odot \tilde{P}_i)$ have multiplicity 1 and
$\tilde{P}_i$ is the projection \lq conjugated\rq\ to $P_i$.

\section{Woronowicz' external algebra for $\mc{O}(O_q(3))$}

Throughout the section we fix a positive integer $k$. Set $\eps=0$  if
$k$ is even and $\eps=1$ if $k$ is  odd. 

\addvspace{\baselineskip}

{\sc Theorem.}
{\em Let $\A$ be the quantum group $O_q(N)$, $\Gamm$ one of the
$N^2$-dimensional bicovariant first order differential calculi $\Gamm_+$ or
$\Gamm_-$ on $\A$ and $q$ be transcendental. Let $\gd$ denote
Woronowicz' external algebra over $\Gamm$. 
\\
Then for each positive integer $k$ there exists a nonzero $k$-form in $\gd$.
In other words $\gd$ is an infinite differential complex.
}

\addvspace{\baselineskip}

We sketch the main steps of the proof.
{}From (\ref{eq-bk}) and (\ref{eq-Tpunkt}) and for $k=3$ we obtain
for example
\begin{displaymath}
\setlength{\unitlength}{1cm}
\raisebox{.4cm}{$\dot{\sigma}_{12}=$\ }
\begin{picture}(1.5,1)
\put(0,1){\vector(1,-1){1}}
\put(1,1){\line(-1,-1){.4}}
\put(.4,.4){\vector(-1,-1){.4}}
\put(1.5,1){\vector(0,-1){1}}
\end{picture}
\raisebox{.4cm}{\ $\odot$\ }
\begin{picture}(1.5,1)
\put(0,0){\vector(1,1){1}}
\put(1,0){\line(-1,1){.4}}
\put(.4,.6){\vector(-1,1){.4}}
\put(1.5,0){\vector(0,1){1}}
\end{picture}
\raisebox{.4cm}{\,,}
\qquad\qquad
\raisebox{.4cm}{$\dot{\sigma}_{23}=$\ }
\begin{picture}(1.5,1)
\put(.5,1){\vector(1,-1){1}}
\put(1.5,1){\line(-1,-1){.4}}
\put(.9,.4){\vector(-1,-1){.4}}
\put(0,1){\vector(0,-1){1}}
\end{picture}
\raisebox{.4cm}{\ $\odot$\ }
\begin{picture}(1.5,1)
\put(.5,0){\vector(1,1){1}}
\put(1.5,0){\line(-1,1){.4}}
\put(.9,.6){\vector(-1,1){.4}}
\put(0,0){\vector(0,1){1}}
\end{picture}
\raisebox{.4cm}{\,.}
\end{displaymath}
Similarly, Woronowicz' formula for $A_k$ gives
$\dota:=\sum_{w\in S_k}(-1)^{\ell(w)}T_w\odot T^\fettc_w$,
where $T_w=\Rda_{i_1}\Rda_{i_2}\cdots \Rda_{i_m}$,
$T_w^\fettc=\Rchm_{i_1}\cdots\Rchm_{i_m}$, and $w=s_{i_1}\cdots
s_{i_m}$ is an expression of the permutation $w$ of length $m$ into a
product of $m$
nearest neighbour transpositions $s_j$. Note that
$\Rchm\in\Mor(u^\fettc\ott u^\fettc)$,
$\Rchm{}^{ab}_{rs}=\Rdam{}^{sr}_{ba}$. 
We construct a nonzero vector $t_k$ in $V^{\ot 2k}$, 
$V=\C^N$, with $\dota t_k=\tau_k t_k$, where $\tau_k$ is a nonzero 
eigenvalue of
$\dota$. Using the decomposition of $u^{\ot k}$ (resp. of $u^{\fettc
\ot k}$) into irreducible
subcorepresentations we get a decomposition of the right coaction of $\A$ on
$\Gamm^{\ot k}$ into  smaller components $\pi_\lambda\odot\pch$.
Here $\lambda$
is a partition of $k-2f$, $f=0,\dots,[k/2]$, and 
$\pi_\lambda$ (resp. $\pch$) stands for the irreducible
subcorepresentation of $u^{\ot k}$ (resp. of $u^{\fettc\ot k}$)
determined by the Young diagram $\lam$ (resp. $\mu$).

{\em Step 1.} Throughout we only consider the case $\lam=(\eps)$ and
$\mu=(k)$. The corresponding minimal central idempotents of $\Mor(u^{\ot k})$
(resp.\ of  $\Mor(u\cont{}^{\ot k})$)
are denoted by $z_{(\eps)}$ (resp.\ by $\zchk$).
Using
$\zchk\Rchm=q^{-1}\zchk$, the
antisymmetrizer $\dota$ reduces to
\begin{displaymath}
\dota(z_{(\eps)}\odot \zchk)=(z_{(\eps)}\odot\zchk)\sum_{w\in
  S_k}(-q)^{-\ell(w)}T_w\odot 1.
\end{displaymath}
Hence it suffices to work with 
$\da:=z_{(\eps)}\sum_{w\in  S_k}(-q)^{-\ell(w)}T_w$ in the simple component
$z_{(\eps)}\Mor(u^{\ot k})$ of $\Mor(u^{\ot k})$. 
Obviously, $\da$ acts on $M_k:=\Mor(u^\eps,u^{\otimes k})$ by composition
(here $u^0=1$ denotes the trivial corepresentation).

{\em Step 2.} We define elements $e_m\in\Mor(u^\eps, u^{\ot k})$
and  $e^m\in\Mor(u^{\ot k},u^\eps)$, $m=1,\ldots,k$
by $e^1=e_1=I$, $e_2=(C^{ab})$,
$e^2=((C^{-1})_{ab})$ and by recursion formulae $e_{k+2}=e_k\ot e_2$
and $e^{k+2}=e^{k+2}\ot e^2$.
Set $\alpha_1=1$, $\alpha_2=1-\RM\QM$, and $r=q^{N-1}$.
Further, we abbreviate
$\alpha_k=\beta_{k,k-1}\alpha_{k-1}$, where 
\mbox{$\beta_{k,k-1}=(1+\cdots+q^{-k+2})\alpha_2$} for $k$ even
and $\alpha_k=\alpha_{k-1}$ for $k$ odd.
Let $t_k=\alpha_k^{-1}\da e_k$.

\addvspace{\baselineskip}

{\sc Lemma.}
{\em {\em (i)} The  antisymmetrizer $\da$ has rank one on the module $M_k$. 
The unique up 
to scalars image is $t_k$. The minimal polynomial of $\da$ on $M_k$ is 
$$\da(\da-\tau_k)=0.$$
\\
{\em (ii)}  We have
$$
(e^2)_{i, i+1}t_k=\gamma_k t_{k-2},\qquad i=1,\dots,k-1,
$$
where \mbox{$\gamma_k=\alpha_2(q-q^{-1})^{-1}(q+q^{-k+1}r)$}
for $k$ odd and $\gamma_k=\gamma_{k-1}$ for $k$ even.
}

\addvspace{\baselineskip}

{}From (ii) it follows that $e^k t_k=\gamma_k\gamma_{k-2}\cdots\gamma_\eps$
is nonzero. Consequently, $t_k$ 
is a nonzero vector. Note that for the quantum group $Sp_q(N)$ this is no 
longer true since $r=-q^{N+1}$. Then we have
$\gamma_N=0$ and moreover $t_N=0$.

{\em Step 3.} It remains to prove that the eigenvalue $\tau_k$ of $\da$ to the 
eigenvector $t_k$ is nonzero: we compute for each $k$
a polynomial $p_k(x,y)$ with integer coefficients 
such that $\tau_k=p_k(\QM,\RM)$, $p_k(0,0)=1$, and $p_k(1,1)=0$. Since $q$ is 
transcendental, $\tau_k\ne0$. The explicit values for small $k$ are
\begin{eqnarray*}
\tau_2&\!\!=&\!\!1-\QM\RM,\quad
\tau_3=1+2q^{-2}-2q^{-1}\RM-q^{-3}\RM,\\
\tau_4&\!\!=&\!\!(1+q^{-2})(1-\QM\RM)\tau_3,\\
\tau_5&\!\!=&\!\!(1+q^{-2})((1+3q^{-2}+6q^{-4}+5q^{-6})\\
&&-\RM(4\QM+11q^{-3}+11q^{-5}+q^{-7}) 
 +r^{-2}(5q^{-2}+6q^{-4}+3q^{-6}+q^{-8})).
\end{eqnarray*}

\section{Further Bicovariant Bimodules for $\mc{O}(O_q(3))$}

Now let $\A=\mc{O}(O_q(3))$ and $v=u$. Further, we suppose that
$q$ is a transcendental complex number.
The matrix of the braiding $\sigma$ with respect to
the basis given above is of the form
$\sigma =b_2(\Rda\odot \Rupm)b_2^{-1}$. It has 7 eigenvalues
(see also \cite{a-CSchWW1}):
$1$, $q^3$, $q^{-3}$, $-q^2$, $-q^{-2}$, $-q$, and $-q^{-1}$.
We have $u\cont \cong u$ and
$\Mor(u^{\otimes k})$ is a factor algebra $B_k$ of the Birman-Wenzl-Murakami
algebra. The algebra $B_2$ has three projections: $P_+$, $P_-$, and $P_0$.
Let us first consider Woronowicz' external algebra.

$a)$ $\mc{S}_{L,1}:=\langle \ker A_k\,|\,k\geq 2 \rangle$.
{}From the previous section we have $\dim \Gamm_{L,1}^\wedge=\infty$.

The following table (which is valid for $\mc{O}(O_q(N))$ and
$\mc{O}(Sp_q(N))$ as well)
gives the projections onto nonzero bicovariant
subbimodules of $\Gamm^{\wedge 3}_{L,1}$ and the corresponding
eigenvalues of $A_3$.\\
\begin{tabular}{ll}
& \\
$\ytt111030000 \odot \ytt300011100$ & $q^{-3}[2][3]$\\
& \\
$\ytt300011100 \odot \ytt111030000$ & $q^3[2][3]$\\
& \\
$\ytt210021000 \odot \ytt210021000$ & $2[3]$\\
& \\
$\ytt300011100 \odot \ytt100010000$ & $1+2q^2-2qr-q^3r$\\
& \\
$\ytt100010000 \odot \ytt300011100$ & $1+2/q^2-2/(qr)-1/(q^3r)$\\
& \\
$\ytt111030000 \odot \ytt100010000$ & $1+2/q^2+2r/q+r/q^3$\\
& \\
$\ytt100010000 \odot \ytt111030000$ & $1+2q^2+2q/r+q^3/r$\\
& \\
$\ytt210021000 \odot \ytt100010000$ & $[3](r+1)-Q(r-1)$, \quad
$-[3](r-1)-Q(r+1)$\\
& \\
$\ytt100010000 \odot \ytt210021000$ & $-[3](1/r-1)+Q(1/r+1)$, \quad
$[3](1/r+1)+Q(1/r-1)$\\
& \\
$\ytt100010000 \odot \ytt100010000$ &
$\lambda_1$, $\lambda_2$, $\lambda_3$, $\lambda_4$.\\
&
\end{tabular}\\
Here we used the $q$-numbers $[2]=q+1/q$, $[3]=q^2+1+q^{-2}$ and
the abbreviations $Q=q-1/q$, $r=q^{N-1}$ for $\mc{O}(O_q(N))$ and
$r=-q^{N+1}$ for $\mc{O}(Sp_q(N))$.
The complex numbers $\lambda_i$, $i=1,2,3,4$ are zeros of the equations
\begin{eqnarray*}
\lambda_{1,2}^2-2((r-1/r)^2+Q(Q^2+1)(r-1/r)+2Q^2)&\!\!=&\!\!0,\\
\lambda _{3,4}^2-2([3]-Q(r-1/r))\lambda _{3,4}-2(r-1/r)((r-1/r)+Q[3])
&\!\!=&\!\!0.
\end{eqnarray*}
In particular if we only require that $q$ is not a root of unity
the vector spaces $\mc{S}^3_{L,1}$ may have different dimensions
for different $q$'s.

$b)$ $\mc{S}_{L,2}:=\langle \ker(I-\sigma) \rangle$.

With our formulation it means
\begin{displaymath}
\mc{S}_{L,2}=\langle b_k(P_+\odot P_+ +P_-\odot P_- +P_0\odot P_0)b_k^{-1}
(\Gamm_L^{\otimes 2})\rangle.
\end{displaymath}

The dimensions of
$\Gamm_{L,2}^{\wedge k}=\Gamm_L^{\otimes k}/{\mc{S}_{L,2}^k}$, $k\leq 6$,
are
\begin{displaymath}
\begin{array}{l|ccccccc}
k & 0 & 1 & 2 & 3 & 4 & 5 & 6 \\
\hline
\dim & 1 & 9 & 46 & 183 & 628 & 1938 & 5514
\end{array}
\end{displaymath}
The computations were carried out with the help of the computer
algebra program FELIX \cite{a-felix}.
Using the table in the case $a)$ it is easy to see that
$\dim \Gamm^{\wedge k}_{L,2}=\dim \Gamm^{\wedge k}_{L,1}$ for $k\leq 3$.
Obviously, $\dim \Gamm^\wedge_{L,2}\geq \dim \Gamm^\wedge_{L,1}=\infty $
by the previous section.

Let us give some motivation for the next definition.

The bicovariant bimodules for a coquasitriangular Hopf algebra
described above admit a second braiding:
\begin{displaymath}
\tilde{\sigma}=b_2(\Rda\odot \Rup)b_2^{-1},
\end{displaymath}
where $\Rda{}^{ij}_{kl}=\br (u^j_k,u^i_l)$ and
$\Rup{}^{ij}_{kl}=\br (u\cont{}^j_k,u\cont{}^i_l)$.
For $\A=\mc{O}(SL_q(N))$ the matrix $\tilde{\sigma}$ has eigenvalues
$-1$, $q^2$, and $q^{-2}$ and
\begin{displaymath}
\ker(I-\sigma)=\im(I+\tilde{\sigma}).
\end{displaymath}
For $\A=\mc{O}(O_q(3))$ the operator $\tilde{\sigma}$ commutes with $\sigma$
and its eigenvalues are
$q^2$, $q^{-2}$, $q^{-4}$, $q^{-1}$, $-1$, $-q^{-3}$.

$c)$ $\mc{S}_{L,3}:=\langle \im(I+\tilde{\sigma})\rangle$.
This definition is equivalent to
$\mc{S}_{L,3}=\langle \ker(I-\sigma),\ker(q^3-\sigma),\ker(q^{-3}-\sigma),
\ker(q+\sigma),\ker(1/q-\sigma)\rangle$.
The dimension of the vector space
$\Gamm_{L,3}^{\wedge k}=\Gamm_L^{\otimes k}/{\mc{S}_{L,3}^k}$ is given by
the following table:
\begin{displaymath}
\begin{array}{l|ccccc}
k & 0 & 1 & 2 & 3 & 4 \\
\hline
\dim & 1 & 9 & 30 & 39 & 0
\end{array}
\end{displaymath}

The corresponding projections are

\begin{tabular}{ll}
& \\
$k=0:$ & $\emptyset \odot \emptyset $\\
& \\
$k=1:$ & $\ytt100010000 \odot \ytt100010000$\\
& \\
$k=2:$ & $\ytt200011000 \odot \ytt110020000 +
\ytt110020000 \odot \ytt200011000$\\
& \\
$k=3:$ & $\ytt300011100 \odot \ytt111030000 +
\ytt210021000 \odot \ytt210021000 +
\ytt111030000 \odot \ytt300011100$\\
& 
\end{tabular}

It was already suggested by Carow-Watamura et al.\ \cite{a-CSchWW1} to prefer
the following choice:

$d)$ $\mc{S}_{L,4}:=\langle \ker(I-\sigma)\oplus \ker(q^3-\sigma)
\oplus \ker(q^{-3}-\sigma)\rangle$.

Then the dimension of
$\Gamm_{L,4}^{\wedge k}
=\Gamm_L^{\otimes k}/{\mc{S}_{L,4}^k}$
becomes
\begin{displaymath}
\begin{array}{l|cccccc}
k & 0 & 1 & 2 & 3 & 4 & 5\\
\hline
\dim & 1 & 9 & 36 & 54 & 1 & 0
\end{array}
\end{displaymath}
and the corresponding projections are

\begin{tabular}{ll}
& \\
$k=0:$ & $\emptyset \odot \emptyset $\\
& \\
$k=1:$ & $\ytt100010000 \odot \ytt100010000$\\
& \\
$k=2:$ & $\ytt200011000 \odot \ytt110020000
+ \ytt110020000 \odot \ytt200011000
+ \ytt110020000 \odot \emptyset
+ \emptyset \odot \ytt110020000$\\
& \\
$k=3:$ & $\ytt300011100 \odot \ytt111030000
+ \ytt210021000 \odot \ytt210021000
+ \ytt111030000 \odot \ytt300011100
+ \ytt111030000 \odot \ytt100010000
+ \ytt100010000 \odot \ytt111030000
+ \ytt100010000 \odot \ytt100010000$\\
& \\
$k=4:$ & $\emptyset \odot \emptyset$\\
&
\end{tabular}

The radical $\mc{R}$ of $\Gamm_{L,4}^\wedge$,
\begin{displaymath}
\mc{R}=\{\rho\in\Gamm_{L,4}^\wedge\,|\,\rho\wedge\theta_{ij}=0
\mbox{ for all }i,j=1,2,3\}
\end{displaymath}
is non-trivial: we have
$\mc{R}\subset\Gamm_{L,4}^{\wedge 3}$ and
$\dim\mc{R}=45$. Hence we can define
$\tilde{\Gamm}_{L}:=\Gamm_{L,4}^\wedge/\mc{R}$
and obtain the dimensions
\begin{displaymath}
\begin{array}{l|cccccc}
k & 0 & 1 & 2 & 3 & 4 & 5\\
\hline
\dim & 1 & 9 & 36 & 9 & 1 & 0
\end{array}
\end{displaymath}

Recall that for the $4D_{\pm}$-calculi on $\mc{O}(SL_q(2))$
the left-invariant external algebra due to Woronowicz is finite dimensional
and has a unique left-invariant form of maximal degree 4 as well.
On the other hand, $\mc{O}(O_q(3))$ is isomorphic to a subalgebra
of $\mc{O}(SL_q(2))$ and the fundamental bicovariant bimodule on
$\mc{O}(O_q(3))$
examined in this paper can be obtained from the 9-dimensional
one on $\mc{O}(SL_q(2))$ (determined by the 3-dimensional
corepresentation of the latter, see \cite{b-KS}).
Moreover, the corresponding exterior algebras have the same dimensions,
since the bicovariant comodules of left-invariant symmetric
$k$-forms are isomorphic for all $k$.
Is this the reason, why the volume form has degree 4?
What happens for $\A=\mc{O}(SL_q(2))$ and the other
bicovariant bimodules (differential calculi)?
Is there a unique form of maximal degree for other quantum groups
such as $\mc{O}(O_q(N))$ or $\mc{O}(Sp_q(N))$, $N\geq 4$?

Let us conclude with two conjectures: firstly, Woronowicz' left-invariant
external algebra is infinite dimensional for the $N^2$-dimensional
bicovariant differential calculi on $\mc{O}(Sp_q(N))$.
Secondly, let $\Gamm$ be an $N^2$-dimensional bicovariant differential
calculus on $\mc{O}(O_q(N))$ or $\mc{O}(Sp_q(N))$, $N\geq 3$. Then
for transcendental values of $q$ the two-sided ideals
$\langle\ker(I-\sigma)\rangle$ and
$\langle\ker A_k\,|\,k\geq 2\rangle$ of the algebra $\Gamm^\otimes$
coincide.

\end{document}